\documentclass[11pt,reqno]{amsart}

\usepackage{amsmath, amsfonts, amssymb, amscd, enumerate}
\newtheorem{theo}{Theorem}[section]
\newtheorem{lem}[theo]{Lemma}
\newtheorem{prop}[theo]{Proposition}

\newtheorem{rem}[theo]{Remark}

\newtheorem{definition}[theo]{Definition}
\newenvironment{pf}{\noindent{\it Proof. }}{$\square$\par\medskip}

\newtheorem{proposition}[theo]{Proposition}

\newcommand{\C}{{\mathbb C}}
\newcommand{\R}{{\mathbb R}}
\newcommand{\Z}{{\mathbb Z}}

\def\CC{{\mathcal C}}

\newcommand{\Jst}{J_{\operatorname{st}}}
\newcommand{\st}{{\operatorname{st}}}
\newcommand{\cc}{\text{c. c.\ }}
\newcommand{\J}{{\mathbb J}}

\newcommand{\D}{{\mathcal D}}
\newcommand{\NC}{{\mathcal N}}
\renewcommand{\L}{{\mathcal L}}
\newcommand{\BL}{{\mathbb L}}
\newcommand{\U}{{\mathcal U}}
\newcommand{\V}{{\mathcal V}}
\newcommand{\F}{{\mathcal F}}
\newcommand{\M}{{\mathcal M}}
\renewcommand{\H}{{\mathcal H}}
\newcommand{\K}{{\mathcal K}}
\newcommand{\I}{{\mathcal I}}

\renewcommand{\=}{\overset{\text{def}}{=}}

\newcommand{\re}{\operatorname{Re}}
\newcommand{\im}{\operatorname{Im}}

\renewcommand{\span}{\operatorname{span}}

\def\sideremark#1{\ifvmode\leavevmode\fi\vadjust{
\vbox to0pt{\hbox to 0pt{\hskip\hsize\hskip1em
\vbox{\hsize3cm\tiny\raggedright\pretolerance10000
\noindent #1\hfill}\hss}\vbox to8pt{\vfil}\vss}}}

\title[An existence theorem for  stationary discs]
{An existence theorem for stationary discs\\ in  almost complex manifolds}

\author{A.Spiro and A.Sukhov}


\begin{document}

\begin{abstract}
An existence   theorem for 
stationary discs of  strongly pseudo-convex domains 
in almost complex manifolds is proved.  More precisely, 
it is shown that, for all points of a suitable neighborhood of the 
boundary and for any 
vector belonging to certain open subsets of the tangent spaces, 
there exists a unique   stationary disc passing through that point and 
tangent to the given vector.
\end{abstract}

\maketitle

\null \vspace*{-.25in}

\section*{Introduction}
Analysis on almost complex manifolds recently became a rapidly growing
research area in modern geometric analysis, due to  the impulse given 
by  the fundamental paper of Gromov \cite{Gr}, where  deep connections between almost
complex and symplectic structures have been discovered. One of the main tools in Gromov's
approach to the study of almost complex manifold  
is given by the {\it pseudo-holomorphic discs\/},  namely
the holomorphic maps from the unit disc $\Delta \subset \C$ 
into the considered almost complex manifold $(M, J)$. When the almost complex structure 
$J$ is integrable, such pseudo-holomorphic discs  coincide with  the well known
 holomorphic discs of a complex manifold and they represent 
one of the central objects of complex analysis in several variables. The theory of classical 
holomorphic discs is by now well developed, while for what concerns 
the pseudo-holomorphic disc in a generic almost complex manifold, several quite natural 
questions are still open. \par
In the classical case,  the theory of  holomorphic discs is 
tightly related with the studies on the  Kobayashi-Royden metric of complex domains. 
Such metric 
can be easily defined  also for domains in  a general almost complex manifold, but,
up to now,
relatively few thing are known on its properties in this more general setting. Quite recently,
new results 
have been obtained by 
Gaussier and the second author in \cite{GS} and by Ivashkovich and Rosay  in \cite{IR}, concerning
the boundary behavior of the Kobayashi metric on strictly pseudoconvex
domains in almost complex manifolds, a class of domains 
which plays a basic role  in Gromov's theory. In this paper 
we continue those studies focusing on the    almost complex
analogue of Lempert's idea  of {\it stationary discs\/} of  pseudoconvex domains.\par
\medskip
The family of stationary discs consists of a special class of holomorphic discs that are
attached to the boundary of a bounded domain. They were considered 
for the first time in  the celebrated paper 
\cite{Le}.  In that paper, Lempert shows that  the extremal discs for the
Kobayashi  
metric of a strongly convex domain  $D\subset \C^n$, coincide
with the   stationary discs, i.e.  with the holomorphic discs  admitting a meromorphic
lift to the cotangent bundle $T^*\C^n$,  with the boundary attached to
the conormal bundle of  $\partial D$ and with exactly one pole of order one at the origin.
  Lempert's idea of stationary discs
turned out to be a very important and fruitful tool  in a variety of topics
of research, like e.g.   Kobayashi metric,  
 boundary regularity problems, moduli
spaces of bounded domains,   
complex Monge-Amp\`ere equations.  \par
Many authors developed the theory of stationary discs. We limit ourselves to mention 
the  recent  Tumanov's paper \cite{Tu}, where a
very simple definition of stationary discs  (equivalent to Lempert's definition) is given and 
is immediately generalizable
to the almost complex setting.  The wide range of applicability of the
stationary discs, even for questions non related 
to the boundary of pseudoconvex domains,  is made evident  f.i. in \cite{Tu} 
and it represents 
one of the main motivations for  studying them  also  in the almost complex case. 
Other good reasons to develop the theory of almost complex stationary discs
can be extrapolated    from the recent results by Coupet, Gaussier and the second 
author in \cite{CGS}, where the
stationary discs of small deformations 
of the standard complex structure of $\C^n$ are used to  
solve  extension problems  for biholomorphisms
 between  domains  in almost complex manifolds.\par
In the
present paper, we are interested in the existence and uniqueness
problem for  stationary  discs of a 
strictly pseudoconvex domain in an almost complex manifolds $(M,J)$. No 
additional restrictions on the almost complex structure $J$ are assumed. 
Our main result is the following (for the definition of $\Gamma$-stationary see
 Def. \ref{stationarydiscs} below).\par
\medskip
\begin{theo} \label{mainresult} Let $D$ be a smoothly bounded   strongly 
pseudoconvex domain   in an almost complex manifold $(M, J)$. 
Then there exists a neighborhood $\U$ of the boundary $\Gamma = \partial D$ such that for any $q\in \U \cap D$ 
and any vector $v$ of a suitable open neighborhood of a codimension one
complex subspace $V_q \subset T_q M$, there exists a unique $\Gamma$-stationary 
$J$-holomorphic  disc $f: \bar \Delta \to \U \cap \bar D$, with 
$f(0) = q$ and $df\left(\left.\frac{\partial}{\partial x}\right|_0\right) \in \R v$.
\end{theo}
The proof  gives also precise information on the complex subspaces
$V_q \subset T_q M$, $q\in \U\cap D$.  In fact, let us fix a point $p \in \Gamma = \partial D$
and let $[q_0,p]$ be a segment which endpoints $p$ and $q_0 \in D$ which is orthogonal
to $\partial D$ with respect to some Riemannian metric on $M$.Then there exists a neighborhood $\U$ of 
$\partial  D$ which satisfies the claim of theorem, and a neighborhood 
$\U' \subset U$ 
of $p$ on which:
\begin{itemize}
\item[(a)] there is a system of coordinates centered at $p$;
\item[(b)]  $\Gamma$ and $J$ have coordinate expressions in the so-called ``standard form" (for
the definition, see  \S 2).  
\item[(c)] $[q_0,p] \cap \U'$ is a part of the coordinate axes $z^1 = ...= z^{n-1} = Im z^n = 0$.
\end{itemize}

Under the identification of 
 $\U'$ with an open subset of $\C^n$ determined by such coordinates, 
 for any $q \in [q_0,p] \cap U'$ {\it the subspaces $V_q$
 are  given by the complex $(n-1)$-dimensional spaces 
 parallel (in $\C^n$) to the holomorphic tangent space of $\partial D$ 
at $p = 0$\/}.  Of course, a different choice of the considered  Riemannian metric 
corresponds to different possibilities for subspaces $V_q$.\par
 \medskip
We   point out that, in the standard case,  the results of Lempert  (\cite{Le})
and M.-Y. Pang (\cite{Pa}) give the existence 
of stationary discs,  
attached to ``convexifiable"  regions of $\partial D$. This gives immediately
the existence of stationary discs attached to sufficiently small
open sets  of the boundary of a strongly pseudoconvex 
domain. Our result can be consider as an extension of this  fact
in the almost complex setting. \par
\medskip
The main ingredient  of our proof consists in a scaling argument, 
which allows to consider  the Riemann-Hilbert problem,  defining the  small stationary 
discs, as a continuous deformation of the Riemann-Hilbert problem which determines
 the  stationary  discs of a suitable  ``osculating"  almost complex structure. 
After proving the existence of stationary discs for the osculating structure, the
result follows as a direct application of the Implicit Function Theorem.\par
\smallskip
Finally, we would like to stress the fact  that our main result
is essentially an application  of  the properties of the local
geometry of strongly pseudoconvex hypersurfaces in the 
``osculating" almost complex structures.
 For  a deeper analysis  of the geometry of such crucial models, 
 we refer the reader to  the recent paper \cite{GS1}.\par
\medskip
The structure of the paper is as follows: after a preliminary section, in \S 2  we give 
the concept of osculating structures and prove the existence of a continuos 
deformation between the Riemann-Hilbert problem for stationary discs of
 the given almost complex structure $J$ and attached to a given 
 boundary $\partial D$ and  the corresponding problem for the osculating 
structure. In \S 3 we give the explicit expressions for 
 the Riemann-Hilbert  problem  determined by the 
osculating structure and in \S 4 we prove the main theorem by showing that the 
deformation described in \S 2 satisfies the hypothesis of the general Implicit 
Function Theorem.\par
\medskip
Everywhere in this paper, we will denote by $\Jst$ the standard complex structure of 
$\C^n$, given by the multiplication by $i = \sqrt{-1}$, and by 
 $\Delta$  the unit disc in $\C$. \par
 \bigskip
 \noindent{\it Acknowledgment\/}. This work was partially done when the second author
 visited the Universities of Camerino and Florence in Spring 2004. He thanks these institutions for 
 hospitality.\par
 \bigskip
\bigskip
\section{Preliminaries}
\bigskip
\subsection{Almost complex structures}
Let $(M,J)$ be an almost complex manifold, i.e.  a  manifold
$M$ endowed with an almost complex structure $J$. We remind that 
an almost complex structure $J$ is   a smooth tensor field  
of type $(1,1)$, such that  $J_x^2 = - Id$ 
at any point $x\in M$. \par
We recall that if $M$ is a complex manifold (i.e. a manifold endowed 
with an atlas of complex coordinate charts $\xi: \U \subset M \to \C^n$ that overlap
biholomorphically), then  it is naturally 
endowed with an almost complex structure $J$. In fact, it suffices to consider  the 
tensor field $J$ which, in any holomorphic system of 
coordinates $\xi = (z^1, \dots, z^n)$, is defined by
$$J_x\left(\frac{\partial}{\partial z^i}\right)  = \sqrt{-1}\frac{\partial}{\partial z^i}\ .\eqno(1.1)$$
An almost complex structure $J$ is called  {\it (integrable) complex structure\/} 
if it is the almost complex structure defined by means of (1.1) by an atlas on $M$
of complex coordinate charts that overlap biholomorphically.\par
The standard  complex structure of $\C^n$ (i.e. the almost complex 
structure determined by the standard coordinates of $\C^n$)
will be always denoted by  $\Jst$. \par
\medskip
Given two almost complex manifolds $(M, J)$ and $(M',J')$, a map 
$f: M \to M'$ is called  {\it $(J,J')$-holomorphic\/} if $df \circ J_x= J'\circ df_x$ at any $x\in M$.
It is well known that, in case $M$ and $M'$ are complex manifolds, then $f$ is holomorphic
if and only if the expression of $f$ in two systems of   holomorphic coordinates for $M$ and $M'$  satisfies 
the classical Cauchy-Riemann equations.\par
Any $(\Jst, J)$-holomorphic  map $f: \Delta\to M$  from the unit disc to
an almost complex manifold $(M,J)$ is usually called
 {\it $J$-holomorphic disc\/}.\par
\smallskip
It is well known that for any complex manifold $M$, the cotangent bundle $T^*M$ has a natural
structure of complex manifold and the complex structures $J$ and $\J$ of 
$M$ and $T^*M$, respectively, are such that the natural projection $\pi:T^*M \to M$ 
is holomorphic. This property has been generalized by Sato in \cite{Sa} for 
any almost complex manifold $(M,J)$ and its cotangent bundle $T^*M$. 
We  summarize the key points of Sato's result in the next 
proposition. In the following, given a system of coordinates $\xi = (x^1, \dots, x^{2n}):\U \subset M \to \R^{2n}$ on $M$, we call  ``associated system of coordinates on $\pi^{-1}(\U) \subset T^*M$"
the set of coordinates $\hat \xi = (x^1, \dots x^{2n}, p_1, \dots, p_{2n})$,  which associate 
to any $\alpha  = p_i dx^i \in T_x^*M \subset  \pi^{-1}(\U)$ the coordinates $(x^i)$ of the point $x$ and 
the components $(p_i)$ of $\alpha$ with respect to the basis $dx^i$.  Moreover, if $J^i_j:\U \to \R$ are the real 
functions which give the components of a (1.1) tensor $J = J^i_j \frac{\partial}{\partial x^i} \otimes dx^j$, 
 we denote by $J^i_{j,k}$ the partial derivatives 
$ J^i_{j,k} \= \frac{\partial J^i_j}{\partial x^k}$.\par
\medskip
\begin{prop}\cite{Sa, IK} \label{liftedalmostcomplex} For any almost complex manifold $(M,J)$, 
there exists a unique almost complex structure $\J$ on $T^*M$ with has the 
following properties:
\begin{itemize}
\item[i)] the projection $\pi: T^*M \to M$ is $(\J,J)$-holomorphic;
\item[ii)]  for any $(J,J')$-biholomorphism $f: M\to N$ between two almost 
complex manifolds $(M,J)$ and $(N,J')$, the natural lifted map $\hat f: T^*N \to T^*M$
is $(\J',\J)$-holomorphic;
\item[iii)] if $J$ is a complex structure, then $\J$ coincides with 
the natural complex structure of $T^*M$;
\item[iv)] for  any  system 
of coordinates on the cotangent bundle, which is associated with some
coordinates $\xi = (x^1, \dots, x^{2n}):\U \subset M \to \R^{2n}$ on $M$, if $J^i_j $ denote the components 
of $J$ in the coordinate basis, for any  $\alpha \in T^xM^*$, $x\in \U$,
$$\J|_{\alpha} = J_{i}^a(x) dx^i \otimes \frac{\partial}{\partial x^a} + 
J_{i}^a(x)  d p_a \otimes \frac{\partial}{\partial p_i} +  
 $$
$$ + 
\frac 1 2  p_a\left( - J^a_{i,j}  + J^a_{j,i} + J^a_\ell \left(J^\ell_{i,m} J^m_j - J^\ell_{j,m} J^m_i \right)\right)
d x^i\otimes \frac{\partial}{ \partial p_j}\ .\eqno(1.2)$$
\end{itemize}
\end{prop}
\medskip
The almost complex structure   $\J$  will be called the {\it  canonical 
lift of $J$ on $T^*M$\/}.\par 
\medskip
\subsection{Hypersurfaces and  holomorphic distributions} 
Let   $S \subset M$ be a submanifold of an almost complex manifold $(M,J)$. The
{\it J-invariant} (or {\it $J$-holomorphic\/}) {\it distribution of $S$\/}  
is the collection of  subspaces $\D_x\subset T_xS$, $x\in S$, 
defined by
$$\D_x = \{\ v\in T_x S\ :\ J(v) \in T_xS\ \}\ .\eqno(1.3)$$ 
Notice that the real subspaces $\D_x$ are given by the real parts 
of the vectors $V\in T^\C_x S \subset T_x^\C M$ such that 
$J V = i V$.  If we consider the subbundles $T^{(1,0)} M$ and $T^{(0,1)} M$
of the complexified tangent bundle $T^\C M$, given by the
$(+i)$- and 
$(-i)$-eigenspaces of $J$ in the complexified  tangent spaces $T^\C_xM$, $x \in M$,
we may also say that $\D_x$ is given by the real parts of the vectors in the 
subspaces 
$$\D^\C_x = T^\C_x S \cap T^{(1,0)}_xM\ .$$  
\par
\medskip
Assume now that $\Gamma$ is an  hypersurface in $M$. 
Notice that, in case $\Gamma$ is (locally) defined as the zero set of a smooth real valued function $\rho$
( i.e. 
$\Gamma = \{\ x\in M\ :\ \rho(x) = 0\ \}$), then the 1-form on $M$ defined by 
$$\vartheta_x = \left(d\rho\circ J\right)|_{T_p\Gamma}\ ,$$
is so that
$$\ker \vartheta|_x = \D_x\ .$$
for any $x\in M$. 
We will call any such 
form a {\it defining 1-form for $\D$\/} and  we will call {\it Levi form at $x\in \Gamma$ 
associated with $\vartheta$\/} the quadratic form  
$$\L_x(v) \= - d\vartheta_x(v, Jv)\ ,\qquad \text{for any}\ v\in \D_x \subset T_x \Gamma\ .$$
 Notice that, for any vector field $X$ with values in $\D$
such that $X_x = v$, we have that 
$$\L_x(v) = - d\vartheta_x(X, JX) = - X(\vartheta(JX))|_x +  JX(\vartheta(X))|_x + \vartheta_x([X, JX]) = $$
$$ = \vartheta_x([X, JX])\ , \eqno(1.4)$$
because $\vartheta(X) \equiv \vartheta(J X) \equiv 0$.  \par
\medskip
Notice also that, 
up to multiplication 
by a non-zero real number,
$\L_x$ does not depend on $\vartheta$. Moreover, by  polarization formula, 
we may say that $\L_x$ is the quadratic form associated with the symmetric bilinear form
$(\BL)^s_x$, where $\BL_x$ is the bilinear form defined by 
$$\BL_x: \D_x\times \D_x \to \R\ ,\qquad \BL_x(v,w) = - d\vartheta_x(v, Jw)$$
and $(\BL)^s_x$ is 
the symmetric part of $\L_x$, i.e.
$$(\BL_x)^s(v,w) = - \frac{1}{2} \left(d\vartheta_x(v, Jw) + d\vartheta_x(w, Jv)\right) \ .$$
If $J$ is an integrable  complex structure, 
the bilinear form $\BL_x$ is symmetric
and  
 it coincides with $(\BL_x)^s$. \par
\medskip
\begin{definition}{\rm We  say that an oriented hypersurface $\Gamma\subset M$ is {\it strongly pseudoconvex\/}
if  the Levi form $\L_x$ is positive definite at every point $x\in \Gamma$, for some $\L_x$ 
 associated to a (local) 
defining function $\rho$ of $\Gamma$ such that $d \rho_x(n) < 0$ for any 
vector  $n \in T_x M$, which is considered as  pointing outwards according to the 
orientation of  $\Gamma$.} 
\end{definition}
\medskip 
In all the following, when $\Gamma$ is oriented, we will consider only  local defining functions satisfying  $d\rho|_\Gamma(n) < 0$ for any outwards pointing vector field
$n$.\par
\medskip
We conclude this section 
by recalling the definition of  {\it conormal bundle $\NC^*(\Gamma)$
of a hypersurface $\Gamma \subset M$\/}.
Consider the submanifold $\pi^{-1}(\Gamma) \subset T^*M$. The conormal bundle $\NC^*(\Gamma)$ is defined
as the 
subbundle of $\pi^{-1}(\Gamma)$ given by 
$$\NC^*(\Gamma) \= \{\ \alpha \in T^*_xM\ ,\ x\in M\ :\ \alpha|_{T_x\Gamma} \equiv 0\ \} \subset 
\pi^{-1}(\Gamma) \subset T^*M\ .$$
\medskip
\subsection{Stationary discs of real hypersurfaces in an almost complex manifold}
Let $(M, J)$ be an almost complex manifold. Notice that, on any tangent space $T_x M$ and
any cotangent space $T^*_x M$, we may consider the following natural action of 
$\C$.  For any $a + i b \in \C$,  $v\in T_x M$ and $\alpha \in T^*_x M$ we set 
$$(a+ i b)\cdot v \= a v + b J(v)\ ,\qquad (a+ i b)\cdot \alpha \= a \alpha + b J^*(\alpha)\ ,\eqno(1.5)$$
where $J^*$ is the  complex structure of $T^*_xM$ defined by $J^*(\alpha)(\cdot)  \= \alpha(J(\cdot))$.
Notice also   that such an action of $\C$ commute with the push-forwards and pull-backs of 
holomorphic maps. More precisely, if $f: M \to N$ 
is a holomorphic map between two almost complex manifolds $(M, J)$ and $(N, J')$, then 
$$f_*(\zeta \cdot v) = \zeta \cdot f_*(v) \ ,\quad   f^*(\zeta \cdot \alpha) = \zeta \cdot f^*(\alpha)\ ,$$
for any $v\in TM$ and any $\alpha \in T^*N$.\par 
We may now introduce the crucial concept of  ``stationary disc" in almost complex manifolds.
\bigskip
\begin{definition} \label{stationarydiscs} {\rm \cite{Tu, CGS}
A continuous map  $f: \bar \Delta \to M$ is called {\it disc attached to a hypersurface $\Gamma \subset M$\/}
if $f|_{\Delta}$ is  $J$-holomorphic and 
$f(\partial \Delta) \subset \Gamma$. \par
An attached disc $f$  
is called {\it stationary 
disc of $\Gamma$\/} (or, shortly, $\Gamma$-stationary)
 if there exists a continuous map $\hat f :  \bar \Delta \to T^*M$ such that:
\begin{itemize}
\item[a)] it projects onto $f$, i.e. $\pi\circ \hat f = f$; 
\item[b)] $\hat f|_\Delta$ is $\J$-holomorphic; 
\item[c)]
$\zeta^{-1}\cdot \hat f(\zeta) \in \NC^*(\Gamma)\setminus \{\text{zero section}\}$ for any 
$\zeta \in \partial \Delta$, where ``\ $ \cdot$\ "  denotes   
the action of $\C$ on $T^*M$  defined in (1.5). 
\end{itemize}
If  $f$ is stationary, 
any map   $\hat f$
satisfying a), b) and c),  is called {\it ($\J$-holomorphic) lift of $f$\/} .
} 
\end{definition}
\noindent {\it Remark.}\  We have to mention that, in the literature on stationary discs, 
what  is usually called ``lift"   is
the meromorphic map $f^*(\zeta) = \zeta^{-1}\cdot \hat f(\zeta)$ and 
not the map $\hat f$ (see e.g. \cite{Le, Tu}).  It goes without saying that there is a natural 1-1 correspondence between 
the two kind of lifts, but, in the following arguments,  our new terminology turns out to be 
slightly more efficient. \par
 \bigskip
We conclude this preliminary section, recalling that  the condition for a map $f:\Delta
\to \U$ to be  $J$-holomorphic is equivalent to the condition
$$J \circ df\left(\frac{\partial}{\partial x}\right) = df\left(\Jst
\left(\frac{\partial}{\partial x}\right)\right)\ .\eqno(1.6)$$
In fact,  the necessity is obvious, while the sufficiency follows from 
the following argument (see e.g. \cite{IR}): if (1.6) holds, by considering  the linear 
maps $df_\zeta$ and $J_{f(\zeta)}$ as extended to  $\C$-linear maps between the 
complexified tangent spaces and observing that any vector $v\in    T_\zeta \C$
can be expressed as  $v = (a+ i b) \partial_x $ for some $a+ i b\in \C$, we have
that for any $v$
$$ J \circ df(v) = 
(a+ i b) \cdot \left(J \circ df\left( \frac{\partial}{\partial x}\right)\right) \overset{\text{(1.6)}}= 
 (a+ i b)\cdot 
 df\left(\Jst 
\left( \frac{\partial}{\partial x}\right)\right) =  df\circ \Jst (v)\ ,$$
which is the $J$-holomorphicity condition.
Replacing $\frac{\partial}{\partial x} = \frac{\partial}{\partial \zeta} + 
\frac{\partial}{\partial \bar \zeta}$, condition (1.6) becomes
$$J \circ df\left(\frac{\partial}{\partial \zeta} + \frac{\partial}{\partial \bar \zeta} \right) - 
i  df\left(\frac{\partial}{\partial \zeta}\right) + i df\left(\frac{\partial}{\partial \bar\zeta} \right) = 0\ .$$
In case $M$ is identified with an open subset of  $\R^{2n}Ê\simeq \C^n$ by means 
of a system of coordinates (see also next  \S 2)
and recalling that the multiplication by $i$ corresponds to the standard complex 
structure of $\C^n$, we obtain the following very useful expression for   the condition of $J$-holomorphicity 
$$(J  + \Jst)_{f(\zeta)}\left(df\left(\frac{\partial }{\partial \bar \zeta}\right) \right)+
(J  - \Jst)_{f(\zeta)}\left(df\left(\frac{\partial }{\partial \zeta}\right) \right)= 0\ .\eqno(1.7)$$
\bigskip
\bigskip
\section{Standard forms  and osculating pairs
for almost complex structures and strongly pseudoconvex 
hypersurfaces}
\bigskip
Let $(M, J)$ be an almost complex manifold of real dimension $2n$ and let
$\U\subset M$ be an open subset, which admits a system of coordinates. 
Notice that any
coordinate map $\xi = (x^1, y^1, \dots, x^{n}, y^n): \U \subset M \to \R^{2n} = \C^n$ gives an identification 
of $\U$   with the open subset $\xi(\U) \subset  \C^n = \{\ (z^1, \dots, z^n)\ ,\ z^i \= x^i + i y^i\ \}$ . In particular, this identification allows to consider on $\U$ the standard complex structure 
$$\underset{(\xi)}\Jst = \frac{\partial}{\partial y^i}Ê\otimes 
dx^i - \frac{\partial}{\partial x^i}Ê\otimes 
dy^i = i \frac{\partial}{\partial z^i}Ê\otimes 
dz^i - i \frac{\partial}{\partial \bar z^i}Ê\otimes 
d\bar z^i\ ,\eqno(2.1)$$
that we will call {\it associated with the system of coordinates $\xi$\/}. 
Of course, one should never forget that the identification with an open subset of $\C^n$ and the consequently defined
complex structure (2.1) on $\U$  {\it does\/} depend on the considered
 coordinates $(x^1, y^1, \dots, x^{n}, y^n)$. However,  when 
 no confusion may occur, we will often omit the indication of the system of coordinates $\xi$, denoting 
the complex structure (2.1) simply by $\Jst$ and by
 identifying  the open subset $\U \subset M$ with a subset of $\C^n$ with respect to any
 regard.\par
 \bigskip
\noindent {\bf Notation}.\ In all the sequel, we will also adopt the following convention on indices: {\it latin letters $i$,$j$ etc. 
 will be use to denote indices running between $1$ and $n$, while greek letters $\alpha$, $\beta$, etc.
 will be always used to indicate indices running between $1$ and $n-1$\/}.\par
 \bigskip
 Let $\Gamma \subset M$ be an oriented strongly pseudoconvex hypersurface
 and $\U$ an open neighborhood of $p \in \Gamma$,  identified with an open subset 
 of $\C^n$ by means of a system of coordinates $\xi = (z^1, \dots, z^n): \U \to \R^{2n} = \C^n$
 so that $\xi(p) = 0$.  It is immediately recognizable that it is always possible to choose the coordinates
 so that, at $p = 0$,  we have that $J|_0 =  \Jst|_0$, 
  the tangent space $T_0\Gamma$ coincides with the hyperplane $x^n = 0$ and that
  the $J$-invariant
  subspace of $\D_0 \subset T_0\Gamma$ coincides with the span of the vectors parallel to $z^n = 0$.  In other 
  words, we may assume that $J$ and $\Gamma$ are of the form
   $$\Gamma\ : \ \rho(z, \bar z) = 0\ ,\ \qquad  \text{with}\ \ $$
   \smallskip
   $$\rho(z, \bar z) = 2 \re(z^n) - \re\left(K_{\alpha\beta} z^\alpha z^\beta
  \right) - H_{\alpha \bar \beta}Êz^\alpha \bar z^\beta + O(|z|^3) \ ,\ \ K_{\alpha\beta} = 
  K_{\beta \alpha}\ ,\ \overline{H_{\alpha\bar \beta}} = H_{\beta \bar \alpha}\ ,\eqno(2.2)$$
  and 
  $$J|_z = \Jst|_z +  L_z + O(|z|^2) \ ,\quad \Jst =  
i \left(\frac{\partial}{\partial z^i}\otimes dz^i  - \frac{\partial}{\partial \bar z^i}\otimes d\bar z^i\right)\ ,
   \eqno(2.3)$$
  where $L_z$ denotes a real tensor field of type (1,1) which depends linearly on the real coordinates 
  $x^i = \re(z^i)$ and $y^i = \im(z^i)$. 
  From the fact that $J^2 \equiv - I$, it follows immediately that $L_z \cdot \Jst = - \Jst \cdot L_z$
  at any $z \in \U$ and hence that  the linear map $L_z: T_z\C^n \to T_z\C^n$ (or, more precisely, 
  its $\C$-linear extension $L_z: T^\C_z\C^n \to T^C_z\C^n$) is of the form 
  $$L_z =
   \left( L^i_{\bar j k}z^k + L^i_{\bar j\bar k} \bar z^k\right)\frac{\partial}{\partial z^i}\otimes d\bar z^j
 + \cc\eqno(2.4)$$
 (here and in all the following we will always write "$\cc$" to indicate the 
 complex conjugate terms of the previous expression).\par
 The following lemma shows that the difference between the Levi forms of $\Gamma$ at $p = 0$,
computed  using $\Jst$ and $J$ is determined exactly by the coefficients $L^n_{\bar j k}$
and $L^n_{\bar j \bar k}$.\par
\medskip
\begin{lem} For any real vector 
$v = v^\alpha \left.\frac{\partial}{\partial z^\alpha}\right|_0 + \overline{v^\alpha}
\left.\frac{\partial}{\partial \bar z^\alpha}\right|_0 \in \D_0 \subset T_0 \C^n$, let us denote by $\L_{\st}(v)$ and $\L(v)$ the value  at $0$ of the Levi form with respect to
$\Jst$ and $J$, respectively,  associated 
with a defining form $\vartheta$ so that $\vartheta|_0 = d\rho \circ J|_0 = d\rho \circ \Jst|_0$. Then
$$\L(v) = \L_{\st}(v) + 2 i \overline{v^\alpha} v^\beta \left(L^n_{\bar \alpha \beta} -
\overline{L^n_{\bar \beta \alpha}} \right)
\ .\eqno(2.5)$$
\end{lem}
\begin{pf}  We want to  compute $\L(v)$ using (1.4). We consider a smooth real
vector field $X_z = X^i(z) \frac{\partial}{\partial z^i}Ê+ 
\overline{X^i(z)} \frac{\partial}{\partial \bar z^i}$, defined  on a neighborhood of the origin and  such that $X|_0 = v$ and so that, 
at any point $z \in \Gamma$, it takes values in the $J$-holomorphic distribution $\D$. Then, 
using (1.4),  
$$\L(v) = - d\vartheta_0(X, JX) = \vartheta_0([X, JX]) = (d\rho\circ \Jst)_0([X,J X]) \ .$$
Notice that the higher order terms in (2.3)  give no contribution to the vector 
$[X,J X]_0 \in T_0\C^n$ and hence that  $[X, JX]_0 = [X, \Jst X]_0 +  [X, L X]_0$. In particular, we may write 
\bigskip
$$\L(v) = \left(d\rho\circ \Jst\right)_0([X, \Jst X]) + (d\rho\circ \Jst)_0([X,L X]) =$$
\smallskip
$$ =  \left(d\rho\circ \Jst\right)_0([X, \Jst X])  + i (dz^n - d \bar z^n)[X,L X]_0\ .\eqno(2.6)$$
\par
\bigskip
\noindent
Now,  let us denote by  $\vartheta_\st$ the 1-form defined on a neighborhood of $0$ as
$\vartheta_\st \= d\rho\circ \Jst$.  Observe  that the restriction of $\vartheta_\st$ to the tangent spaces of $\Gamma$ gives a defining 1-form 
for the $\Jst$-holomorphic distribution $\D_\st$ of $\Gamma$  and that 
$- d\vartheta_\st|_0(v, \Jst v) = \L_\st(v)$.  So,  we may write also that 
\bigskip
$$\left(d\rho \circ \Jst\right)_0([X, \Jst X])  = \vartheta_\st|_0([X, \Jst X]) = $$
\smallskip
$$ = 
 - d\vartheta_{\st}|_0(X, \Jst X) + \left.X\left(\vartheta(\Jst X)\right)\right|_0 - 
\left.\Jst X\left(\vartheta(X)\right)\right|_0 =  $$
\smallskip
$$ = \L_\st(v) -  \left.X\left(X(\rho)\right)\right|_0 - 
 \left.J_{\st}X\left(J_{\st}X(\rho)\right)\right|_0\ .\eqno(2.7)$$
 \par
 \bigskip
 By construction, at all points of $\Gamma = \{ \rho = 0\}$,  the vector field $X$
 belongs to $\D \subset T\Gamma$. From this we get that $X(\rho)|_z = 0$ and 
 $JX (\rho)_z =  \Jst X(\rho)_z + L X(\rho)|_z+ O(|z|^2) = 0$ at any $z\in \Gamma$. 
Now, recalling that  $J_{\st}X|_0 \in T_0\Gamma$,  we have the following values 
for the directional derivatives  $ \left.X\left(X(\rho)\right)\right|_0 $ and 
$ \left.J_{\st}X\left(J_{\st}X(\rho)\right)\right|_0$
\bigskip
 $$\left.X\left(X(\rho)\right)\right|_0 = 0\ ,\qquad 
 \left.J_{\st}X\left(J_{\st}X(\rho)\right)\right|_0 = - \left.J_{\st}X\left(L X(\rho)\right)\right|_0\ .\eqno(2.8)$$
 \par
 \bigskip\noindent
By means of (2.8) and (2.7) we get that  (2.6) is equal to 
\bigskip
$$\L(v) =  \L_\st(v)  +
 \left.\Jst X\left(L X(\rho)\right)\right|_0  + 
 i (dz^n - d \bar z^n)[X,L X]_0   \ .\eqno(2.9)$$
 \par
 \bigskip
 \noindent
On the other hand,  recalling that the  $\C$-valued functions $X^i(z)$ satisfy $X^n(0) = 0$ and 
$X^\alpha(0) = v^\alpha$ , we may easily compute what we need, i.e. 
\bigskip
$$[X,L X]_0 = X_0(L X) - (L_0 X_0)(X) =  \left.v \left(L_z\left(X^i(z) \frac{\partial}{\partial z^i}Ê+ 
\overline{X^i(z)} \frac{\partial}{\partial \bar z^i}\right)\right) \right|_0=  $$
\smallskip
$$ = \left( v^\alpha \overline{v^{\beta}} L^i_{\bar\beta \alpha}   + \overline{v^\alpha} \overline{v^{\beta}} L^i_{\bar\beta \bar \alpha}  \right) \left.
\frac{\partial}{\partial z^i}\right|_0 + 
\left( \overline{v^\alpha} v^{\beta}  \overline{L^i_{\bar\beta \alpha}}   + v^\alpha v^{\beta} \overline{L^i_{\bar\beta \bar \alpha}} \right) \left.
\frac{\partial}{\partial \bar z^i}\right|_0\  ,\eqno(2.10)$$
\bigskip
$$ \left.\Jst X\left(L X(\rho)\right)\right|_0 =  
i v^\alpha \overline{v^\beta} L^n_{\bar \beta \alpha} - i \overline{v^\alpha} \overline{v^\beta}
 L^n_{\bar \beta\bar \alpha} - i \overline{v^\alpha} v^\beta \overline{L^n_{\bar \beta \alpha} }
+ i v^\alpha v^\beta \overline{L^n_{\bar \beta\bar \alpha}}
\  .\eqno(2.11)$$
\par \bigskip
\noindent
Substituting (2.10) and (2.11) into (2.9), formula (2.5) follows.
\end{pf}
 \medskip
 From the previous lemma, it follows that  if the coefficient  $L^n_{\bar \alpha \beta}$ 
 are vanishing, then  the Levi form at $0$  of $\Gamma$ with respect to $J$ coincides 
 with the Levi form of $\Gamma$ with respect to  $\Jst$. This fact suggests  the following definition
 and motivates the next proposition.\par
 \medskip
  \begin{definition} \label{standard form}{\rm  We say that $J$ and 
 $\Gamma$ are  {\it in standard form in the system of coordinates $(z^1, \dots, z^n)$\/} if:
 \begin{itemize}
 \item[a)]  $J$  is of the form 
 $$J = \Jst + \left( L^i_{\bar j k}z^k + L^i_{\bar j\bar k} \bar z^k\right)\frac{\partial}{\partial z^i}\otimes d\bar z^j
 + \overline{ \left(L^i_{\bar j k}z^k + L^i_{\bar j\bar k} \bar z^k\right)}\frac{\partial}{\partial \bar z^i}\otimes dz^j 
 + O(|z|^2)\ ,$$
 with coefficients $L^i_{j\bar k}$ so that  
 $$L^n_{\bar \alpha \beta} = 0\quad \text{and}\quad 
  L^n_{\bar \alpha \bar \beta} = - L^n_{\bar \beta \bar \alpha}\quad  \text{for any}\  1 \leq \alpha, \beta\leq n-1\ ; \eqno(2.12)$$
 \item[b)] $\Gamma$ admits a defining function on a neighborhood of the origin of the form 
 $$\rho(z,\bar z) = 2 \re(z^n) - \sum_{\alpha=1}^{n-1} |z^\alpha|^2  + O(|z|^3)\ .\eqno(2.13)$$
 \end{itemize}
 \ \/}
 \end{definition}
 \medskip
 \begin{prop} \label{proposition standard} For any $p\in \Gamma$, there exists a neighborhood $\U$ of $p$ and a system of complex coordinates $\xi = (z^1, \dots, z^n)$ with $z(p) = 0$,  in which 
 $J$ and $\Gamma$ are in standard form. In particular, up to a positive scalar multiple, 
 the Levi form 
 of $\Gamma$ with respect to $J$ coincides with the Levi form with respect to $\Jst$ at  $p = 0$ 
 (as usual, $\Jst$ denotes 
 the complex structure associated with the coordinates $ (z^1, \dots, z^n)$). 
 \end{prop}
 \begin{pf} There is no loss of generality if we assume that $\Gamma$ and $J$ are
 of the form (2.2) and (2.3).  Now,  consider  the change of coordinates 
 $$z^\alpha =w^\alpha\ ,\qquad z^n = w^n + \frac{i}{2} L^n_{\bar \alpha  \beta}
 \bar w^\alpha w^\beta +  \frac{i}{4} L^n_{\bar \alpha \bar  \beta}
 \bar w^\alpha \bar w^\beta 
 \ .$$
 The defining function $\rho$ remains of the form (2.2) even when it is written 
 in terms of $w$ and $\bar w$, while $J$ becomes 
 of the form 
 $$J|_w =
i \left(\frac{\partial}{\partial w^i}\otimes d w^i  - \frac{\partial}{\partial \bar w^i}\otimes d
\bar w^i\right) +  L_{w} + O(|w|^2) - $$
$$- \left(L^n_{\bar \alpha \beta}  w^\beta + \frac{1}{2} ( L^n_{\bar \alpha \bar  \beta} 
+  L^n_{\bar \beta \bar  \alpha}) \bar w^\alpha \right)
 \frac{\partial}{\partial w^n}\otimes 
d\bar w^\alpha -$$
$$- \left(\overline{L^n_{\bar \alpha \beta}} \bar w^\beta + \frac{1}{2} ( 
\overline{L^n_{\bar \alpha \bar  \beta} }
+  \overline{L^n_{\bar \beta \bar  \alpha}})  w^\alpha \right)
\frac{\partial}{\partial \bar w^n}\otimes d w^\alpha\ .$$
and it satisfies  condition a) of Definition 
\ref{standard form}.  So, assuming that a) holds, from (2.5) we have that $\L_\st|_0 = 
\L|_0$ and hence that the matrix $H_{\alpha \bar \beta}$ is positive definite. 
Notice also that any change of coordinates of the form $z^\alpha = U^\alpha_\beta w^\beta$, 
$z^n = w^n$, 
gives an expression for $J$, which still satisfies a), while it changes the coefficients  $H_{\alpha \bar \beta}$ into 
the coefficients $H'_   {\alpha \bar \beta} = U^\gamma_\alpha H_{\gamma\bar \delta}
\overline{U^\delta_\beta}$. Hence, by means a 
 suitable change of coordinates,  we may always 
 assume that $H_{\alpha \bar \beta}Ê= \delta_{\alpha \beta}$. Finally, by performing 
 the change of coordinates $z^\alpha = w^\alpha$ and $z^n = w^n - \frac{1}{2}K_{\alpha \beta} w^\alpha w^\beta$, 
 which also leaves the property a) unchanged, we obtain a system of 
 coordinates in which also b) of Definition \ref{standard form} is satisfied.
 \end{pf}
\bigskip
We conclude this section introducing the crucial concept of osculating pairs 
and the proof that the existence problem for small stationary disc can be reduced
to the analysis of  deformations of the osculating structures.\par 
\medskip
\begin{definition}{\rm Let $J$ and $\Gamma$ be in standard form in a system of complex 
coordinates $(z^1, \dots, z^n)$.   We call {\it osculating pair for 
$(J, \Gamma)$ at $p=0$\/} the pair $(J^0, \Gamma^0)$ given by the hypersurface
(locally $\Jst$-biholomorphic to the unit sphere)
$$\Gamma^0 \ :\  \rho^0(z,\bar z) = 2 \re(z^n) - \sum_{\alpha=1}^{n-1} |z^\alpha|^2\ ,
\eqno(2.14)$$
and the almost complex structure 
$$J^0 \= \Jst + 
A_{\bar \alpha \bar \beta} \bar z^\beta \frac{\partial}{\partial z^n} \otimes d \bar z^\alpha  + 
\overline{A_{\bar \alpha \bar \beta}} z^\beta \frac{\partial}{\partial \bar z^n} \otimes d z^\alpha
\ ,\quad \text{with}\ A_{\bar \alpha\bar \beta} = L^n_{\bar \alpha \bar \beta}\eqno(2.15)$$
 where we denoted by  $L^i_{\bar j\bar k}$ and $L^i_{\bar j k}$
  the coefficients of the linear part of $J$ as  in (2.4).\par
}
\end{definition}
\begin{rem}{\rm In case $M$ has real dimension $4$,  the indices $\alpha$ and $\beta$ 
may  assume only the value $1$ and,  by   (2.12), 
$A_{\bar 1 \bar 1} = 0$. In other words, when $\dim_\R M = 4$, 
the only possible osculating pair is  $(\Jst, \Gamma^0)$. If this is the case, 
 all claims of next section  and our main result 
 can be considered as consequences  of  the well-known properties of the 
stationary discs of the unit ball with respect to the standard complex structure $\Jst$. 
}
\end{rem}
\bigskip
Consider now the anisotropic dilations
$$\phi_t: \C^n \to \C^n\ ,\qquad \phi_t(z) = (\frac{1}{t} z^1, 
\dots, \frac{1}{\delta} z^{n-1}, \frac{1}{t^2} z^n)\ ,\qquad \ t \in \R$$
and the pairs $(J^t, \Gamma^t)$ of almost complex structures $J^t \= \phi_t{}_*(J)$ and  hypersurfaces 
 $\Gamma^t \= \phi_t(\Gamma)$.\par
A very simple computation shows that,  for any real value of  $t$,  the 
pair $(J^t, \Gamma^t)$ is in standard form and with $(J^0, \Gamma^0)$ as
osculating pair. Furthermore,  the functions
which give the components  of $J^t$ and the defining functions $\rho^t \= 
\frac{1}{t}(\rho \circ \phi^{-1}_t)$ of the oriented hypersurfaces $\Gamma^t$,  tend uniformly 
on compacta to the components of $J^0$ and the defining function $\rho^0$.
Finally, it can be checked that the partial differential equations and the boundary 
conditions which define the lifts of the stationary discs of the hypersurfaces $\Gamma^t$
with respect to the structures $J^t$ depends continuously on $t$ and they can 
be considered as continuous deformations of the corresponding equations 
and boundary conditions for the lifts of stationary discs 
for $(J^0, \Gamma^0)$.  We formalize this fact in the following proposition, 
which will turn out to be  the key point for the proof of our main result.\par
\medskip
\begin{proposition} \label{deformation} For any $p\in \Gamma$ and any $\epsilon >0$,
 there exists a neighborhood $\U$ of $p$,  
 a system of complex coordinates $\xi = (z^1, \dots, z^n)$ with $z(p) = 0$
 and a smooth family of pairs $(J^t, \Gamma^t)$
 of  almost complex structure $J_t$ and real hypersurfaces $\Gamma^t$, 
 $t \in [-1, 1] \subset \R$, 
 such that
 \begin{itemize}
 \item[i)] $(J^t, \Gamma^t)$ are in standard form for any $t$ and, for all of them, 
 the osculating pair at $p = 0$ coincides with $(J^0, \Gamma^0)$;
 \item[ii)] $(J^1, \Gamma^1) = (J, \Gamma)$ and 
  for any $t\neq 0$, there exists a $(J, J^t)$-biholomorphism  $\phi_t: \U \to \U$ 
 such that $\phi_t(\Gamma) = \Gamma^t$. 
 \end{itemize} 
Moreover,  the partial differential boundary problem which defines the lifts
 of stationary discs of $(J,\Gamma)$ in a sufficiently small neighborhood of 
 the origin is equivalent to an arbitrarily small, 
 continuous deformation of the corresponding  problem for the osculating pair
 $(J^0, \Gamma^0)$
\end{proposition}
\bigskip
\bigskip
\section{The equations for the stationary discs of an osculating pair}
\bigskip
In this section, $(J^0, \Gamma^0)$ will denote a fixed osculating pair, i.e.  a given almost 
complex structure $J^0$ on $\C^n$ of the form (2.15) together with the boundary of the 
Siegel domain  $\Gamma^0$,  defined in (2.14) . We  denote by $(x^1, y^1, \dots, x^n, y^n, u_1, v_1, \dots, u_n, v_n)$  the 
system of coordinates on $T^*\C^n$,  associated with  real coordinates 
$x^i = \re(z^i)$ and $y^i = \im(z^i)$ (see the definition in \S 1.1). This means 
that any 1-form $\alpha \in T^*_z\C^n$ is written in terms of such coordinates as 
$\alpha = u_i d x^i|_z + v_i d y^i|_z$. On the other hand, it is quite useful to consider the $u^i$'s 
and $v^i$'s as real and imaginary parts of some complex coordinates of $T^*\C^n = \C^{2n}$. More 
precisely, 
we will consider on $T^*\C^n$  the complex coordinates 
$(z^1,  \dots, z^n, P_1, \dots, P_n)$, 
with $P_i = \frac{1}{2}(u_i - v_i)$, so that any 1-form $\alpha \in T^*_z\C^n$ will be written
 as 
 \medskip
$$\alpha = P_i d z^i + \bar P_i d\bar z^i\ .$$
\par
\medskip
\noindent
Using (1.2) we may compute the components of the lift $\J^0$ in the real coordinates 
$(x^i, y^i,  u_j, v_j)$  and then re-express $\J^0$ using the complex basis 
associated with the complex coordinates $(z^i, P_j)$. Some tedious but straightforward
computations show that $\J^0$ is 
$$\J^0 =  i \left( \frac{\partial}{\partial z^k}\otimes dz^k -  \frac{\partial}{\partial \bar z^k}\otimes d\bar z^k\right)
  +  A_{\bar \alpha \bar \beta} \bar z^\beta \frac{\partial}{\partial z^n}\otimes d\bar z^\alpha +
  \overline{A_{\bar \alpha \bar \beta}}Ê z^\beta \frac{\partial}{\partial \bar z^n}\otimes d z^\alpha
   +$$
  $$ +  i \left( \frac{\partial}{\partial P_k}\otimes d P_k-  \frac{\partial}{\partial   \bar P_k}\otimes d \bar P_k\right) 
   + \overline{A_{\bar \alpha \bar \beta}} z^\beta \frac{\partial}{\partial \bar P_\alpha}\otimes d  P_n +
A_{\bar \alpha \bar \beta} \overline{z^\beta} \frac{\partial}{\partial  P_\alpha}\otimes d \bar P_n
\eqno(3.1)$$
\par
\bigskip
For the interested reader, we give here some more detailed indications on how (3.1) is obtained. As we mentioned, we have  to apply formula (1.2)
by considering as real coordinates $x^i$ and $p_a$  the coordinates $x^i = Re(z^i)$,  $x^{n+i} = y^i = Im(z^i)$ and 
  $p_a = u_a = \frac{1}{2}Re(P_a)$ and $p_{a + n} = v_a = - \frac{1}{2}Im(P_a)$. In order to simplify the notation in the next formulae, in place of the indices ``$j$" and ``$j+n$" we are going to use the indices $j_r $ and $j_i$, respectively,  to indicate in a more expressive way 
  when we refer to 
  a quantity related to a real  part  (``r") or to an  imaginary part  (``i") of a 
 complex quantity. 
 As usual, all latin indices $j, k$ etc. will be considered running between $1$ and $n$, while greek 
 indices $\alpha$, $\beta$,  etc. will run just between $1$ and $n-1$. Capital letters $A$, $B$, etc. 
 will be used for indices that can 
 be both of the form $j_r$ and of the form $j_i$. \par
 Using these conventions and from the definition of the 
 osculating structure $J^0$, we have that the components  $J^A_B$ of $J^0$ are
 \bigskip
 $$J^{\alpha_r}_{\beta_r} = 0\ ,\ J^{\alpha_i}_{\beta_r}  = \delta_\beta^\alpha \ ,
 \ J^{n_r}_{\beta_r}Ê= \re(A_{\bar \beta\bar \gamma}) x^\gamma + \im(A_{\bar \beta\bar \gamma}) y^\gamma\ ,
 \ J^{n_i}_{\beta_r}Ê= -  \re(A_{\bar \beta\bar \gamma}) y^\gamma + \im(A_{\bar \beta\bar \gamma}) x^\gamma$$
 \smallskip
  $$J^{\alpha_r}_{\beta_i} =  - \delta_\beta^\alpha\ ,\ J^{\alpha_i}_{\beta_i}  = 0 \ ,
 \ J^{n_r}_{\beta_i}Ê= - \re(A_{\bar \beta\bar \gamma}) y^\gamma + \im(A_{\bar \beta\bar \gamma}) x^\gamma\ ,$$
 \smallskip
$$ \ J^{n_i}_{\beta_i}Ê= -  \re(A_{\bar \beta\bar \gamma}) x^\gamma  - \im(A_{\bar \beta\bar \gamma}) y^\gamma\ ,$$
\smallskip
 $$J^{j_r}_{n_r} = 0\ ,\ J^{j_i}_{n_r}  = \delta_n^j \ ,
 \ J^{j_r}_{n_i} = - \delta_n^j \ ,\ J^{j_i}_{n_i}  = 0\ .\eqno(3.2)$$ 
 \par\bigskip
 \noindent
 In particular, the only non trivial values of the partial derivatives $J^A_{B,C}$ are
 \bigskip 
 $$J^{n_r}_{\alpha_r, \beta_r} = \re(A_{\bar \alpha \bar \beta})\ ,\ 
 J^{n_r}_{\alpha_r, \beta_i} = \im(A_{\bar \alpha\bar \beta})\ ,\ 
 J^{n_i}_{\alpha_r, \beta_r} = \im(A_{\bar \alpha\bar \beta})\ ,\ 
 J^{n_i}_{\alpha_r, \beta_i} = - \re(A_{\bar \alpha\bar \beta})\ ,$$
 \smallskip
 $$J^{n_r}_{\alpha_i, \beta_r} = \im(A_{\bar \alpha \bar \beta})\ ,\ 
 J^{n_r}_{\alpha_i, \beta_i} = - \re(A_{\bar \alpha\bar \beta})\ ,\ 
 J^{n_i}_{\alpha_i, \beta_r} = - \re(A_{\bar \alpha\bar \beta})\ ,\ 
 J^{n_i}_{\alpha_i, \beta_i} = - \im(A_{\bar \alpha\bar \beta})\ .$$
 \par\bigskip
 \noindent
 Hence, we may compute the coefficients 
 of the second line terms in (1.3)  and we get that 
 \bigskip
 $$ u_n(- J^{n_r}_{\alpha_r, \beta_r} + 
  J^{n_r}_{\beta_r, \alpha_r} - J^{n_i}_{\alpha_r, \beta_i} + 
 J^{n_i}_{\beta_r, \alpha_i})  + v_n(- J^{n_i}_{\alpha_r, \beta_r} + 
  J^{n_i}_{\beta_r, \alpha_r} + J^{n_r}_{\alpha_r, \beta_i} - 
 J^{n_r}_{\beta_r, \alpha_i}) = 0\ .
 $$
 \par \bigskip
 \noindent
Similarly, we have that all other coefficients  of that line are vanishing. On the other hand, 
 using (3.2) and  expressing the terms of the first line of (1.2) using the complex coordinates 
  $z^i$ and $P_j$, one immediately gets (3.1).\par
  \bigskip
  \medskip
Coming back to (3.1), it is clear that the matrix $(\J^i_j)$,  representing  $\J^0$  in the basis $\left(\frac{\partial}{\partial z^i}, \frac{\partial}{\partial \bar z^i}, 
  \frac{\partial}{\partial P_j}, \frac{\partial}{\partial \bar P_i}\right)$, is 
  $$(\J^i_j) =
   \left(
  \begin{matrix}
   i I_n & {\mathcal A}(\bar z) & 0 & 0\\
  \bar  {\mathcal A}(z) & - i I_n & 0 & 0\\
   0 &  0 &   i I_n & \bar  {\mathcal A}^t(z)\\
   0 & 0 &   {\mathcal A}^t(\bar z) & - i I_n
  \end{matrix}
  \right)\eqno(3.3)
  $$
with 
  $$  {\mathcal A}^{\alpha}_j(\bar z) =  {\mathcal A}^n_n(\bar z)  = 0
  \ , \qquad   {\mathcal A}^n_{\alpha}(\bar z) = A_{\bar \alpha\bar \beta} \bar z^\beta\ .$$
  Recall also that, by (2.12), the coefficients  $A_{\bar \alpha \bar \beta}$ satisfy 
  $$A_{\bar \alpha \bar \beta} = - A_{\bar \beta \bar \alpha}\ \qquad \Rightarrow\qquad A_{\bar \alpha
  \bar \alpha} = 0\quad \text{for any}\ 1 \leq \alpha \leq n-1\ .$$
From (3.3), a direct computation gives the matrix representing the linear map 
$(\J^0 + \Jst)^{-1} \cdot (\J^0 - \Jst)$ in the previous  complex basis. This matrix  is 
\bigskip
$$ \left(
  \begin{matrix}
  0 & - \frac{i}{2}{\mathcal A}(\bar z) & 0 & 0\\
  \frac{i}{2} \bar  {\mathcal A}(z) & 0& 0 & 0\\
   0 &  0 &  0 & \frac{i}{2}\bar  {\mathcal A}^t(z)\\
   0 & 0 &  - \frac{i}{2}{\mathcal A}^t(\bar z) & 0
  \end{matrix}
  \right)$$
  \par
  \bigskip \noindent
  and hence, from (1.7), we obtain that a map $\hat f = (f^1, \dots, f^n, g_1, \dots, g_n): \Delta \to T^*\C^n$ is $\J^0$-holomorphic if and only if it satisfies the following p.d.e. system:
  \bigskip
 $$\frac{\partial f^\alpha}{\partial \bar \zeta} = 0\ ,\qquad
  \frac{\partial f^n}{\partial \bar \zeta}  - \frac{i}{2} A_{\bar\alpha\bar \beta} \overline{f^\beta}
   \overline{\left(\frac{\partial f^\alpha}{\partial \zeta}\right)} = 0\ ,$$
   \medskip
  $$\frac{\partial g_\alpha}{\partial \bar \zeta}  + \frac{i}{2} \overline{A_{\bar \alpha\bar \beta} }f^\beta
   \overline{\left(\frac{\partial g^n}{\partial \zeta}\right)} = 0\ ,\qquad 
   \frac{\partial g_n}{\partial \bar \zeta} = 0\ .\eqno(3.4)
  $$
  \par \bigskip \noindent
We now want  to write down the boundary condition for a $\J^0$-holomorphic disc $\hat f$ in order
to be the lift of a stationary disc. Observe that, by (1.5), the action of a complex number $\zeta$ 
on a real  form $\alpha = P_j d z^j + \bar P_j d\bar z^j \in T^* \C^n$ is 
\bigskip
$$\zeta\cdot (P_j dz^ j +\bar P_j d\bar z^ j) \= 
P_j ( \re(\zeta) dz^ j  + \im(\zeta) (dz^ j \circ  J^0)) + \overline{P_j} ( \re(\zeta)d\overline z^ j +
 \im(\zeta) (d\bar z^ j \circ  J^0)) = $$
 \smallskip
$$= \zeta P_j
 dz^ j + \overline{\left(\zeta P_j\right)}
 d\bar z^ j +
 \left(\im(\zeta)  A_{\bar \alpha\bar \beta}  \overline{z^\beta} P_n \right) d\overline{z^\alpha}
+  \left(\im(\zeta)  \overline{A_{\bar \alpha\bar \beta} }z^\beta \overline{P_n} \right) dz^\alpha\ .$$
\par
\bigskip
\noindent
Using the fact that the fibers of $\NC(\Gamma^0)$ are generated over $\R$ by the 1-forms 
$$d\rho^0|_z = d z^n - \sum_{\alpha = 1}^{n-1} \bar z^\alpha d z^\alpha 
+ d\bar z^n - \sum_{\alpha = 1}^{n-1} z^\alpha d\bar z^\alpha\ ,\quad z\in \Gamma^0\ , $$
it is quite immediate to realize that  a $\J^0$-holomorphic disc $\hat f = (f^1, \dots, f^n, g_1, \dots, g_n)$ 
satisfies condition (c) of Definition \ref{stationarydiscs} 
if and only if  for any $\zeta \in \partial \Delta$
there is a real number $\lambda_\zeta \neq 0$ such that
\bigskip
$$2 \re(f^n(\zeta)) - \sum_{\alpha = 1}^{n-1} |f^\alpha(\zeta)|^2 = 0\ ,\eqno (3.5)$$
\bigskip
$$g_i(\zeta) dz^i + \overline{g_i(\zeta)} d\bar z^i = \lambda_\zeta \left(\zeta \cdot d\rho^0|_{f(\zeta)}\right) = $$
\smallskip
$$ = \lambda_\zeta \left(
- \zeta \overline{f^\alpha(\zeta)} dz^\alpha + \zeta d z^n
- \frac{i}{2} (\zeta - \bar \zeta) \overline{A_{\bar \alpha\bar \beta} }f^\beta(\zeta) d z^\alpha -
\phantom{aaaaaaaaaaaaaa} \right.$$
\smallskip
$$\left. \phantom{aaaaaaaaaaaaaa}
- \bar \zeta f^\alpha(\zeta) d\bar z^\alpha + \bar \zeta d \bar z^n-  \frac{i}{2} (\zeta - \bar \zeta)  A_{\bar \alpha\bar \beta} \overline{f^\beta(\zeta)} d\bar z^\alpha 
\right)\ .\eqno(3.6)$$
\par \bigskip
\noindent
This immediately implies that $\lambda_\zeta = \bar \zeta g_n(\zeta) = \zeta \overline{g_n(\zeta)}$
and hence that (3.6) is equivalent to 
\bigskip
$$g_\alpha(\zeta) +  \left(\overline{f^\alpha(\zeta)}
+ \frac{i}{2} \overline{A_{\bar \alpha\bar \beta} }f^\beta(\zeta) \right)g^n(\zeta) -
 \frac{i}{2}  \overline{A_{\bar \alpha\bar \beta} }f^\beta(\zeta)\overline{g^n(\zeta)} = 0\ ,
 \  \bar \zeta g_n(\zeta) - \zeta \overline{g_n(\zeta)} = 0\ .\eqno(3.7)$$
 \par \bigskip \noindent
for any $\zeta \in \partial \Delta$.\par
In brief, we have proved the following lemma.\par
\medskip
\begin{lem} A map $\hat f = (f^1, \dots, f^n, g_1, \dots, g_n): \bar \Delta \to \C^{2n}$
represents the lift of a stationary disc of $\Gamma^0$ if and only if it satisfies 
the generalized Riemann-Hilbert problem given by the p.d.e. system (3.4) and 
the boundary conditions (3.5) and (3.7).
\end{lem}
\medskip
Now, recalling that $A_{\bar 1\bar 1} = 0$, 
by a direct inspection it is possible to check that, for any $a \in \C^*$ and  $\lambda \in R^*$,
the map $\hat f_{a, \lambda}  = (f_a(\zeta); g_{a, \lambda}): \bar \Delta \to \C^{2n}$, 
defined by
\bigskip
$$f_{a}(\zeta) = (a \zeta, 0, \dots, 0, \frac{|a|^2}{2})\ ,\eqno(3.8)$$
\bigskip
$$g_{a, \lambda}(\zeta) = (
- \lambda \bar a, \frac{i \lambda }{2} \overline{A_{\bar 2\bar 1}} a(- \zeta^2 - |\zeta|^2 + 2), 
\dots, \frac{i \lambda }{2} \overline{A_{\overline{n-1}\bar 1}}a(- \zeta^2 - |\zeta|^2 + 2), 
\lambda \zeta)\eqno(3.9)$$
\par \bigskip
\noindent 
is  a lift of the stationary discs $f_a: \bar \Delta \to \C^n$ of $\Gamma^0$. \par
Furthermore, notice  that any map of the form 
$z^n = \tilde z^n$, $z^\alpha = U^\alpha_\beta \tilde z^\beta$, with $U^\alpha_\beta \in \operatorname{U}_n$, leaves $\Gamma^0$ invariant and  sends  $J^0$  into an almost
complex structure $J'{}^{0}$, which is still of the form (3.1). It follows immediately 
that the discs (3.8) and (3.9), computed with the coefficients $A'_{\bar \alpha, \bar \beta}$
of $J'$, are images under the previous transformation
of lifts of stationary discs for $(J^0, \Gamma^0)$.  In particular, 
it follows that for any point 
$z_o = (0, \dots, 0,  z_o^{n})$, with $\rho^0(z_o) > 0$,  and 
any vector $v_o\in \span_\C\left\{ \frac{\partial}{\partial z^\alpha}\ \right\}_{1 \leq \alpha \leq n-1}$, there exists a  
stationary disc $f$ for $(J^0, \Gamma^0)$ so that $f(0) = z_o$ and 
 $f_*\left(\partial_{Re \zeta}|_0\right) = \lambda v_o$ for some $\lambda \in \R^*$. 
 \par
\bigskip
\bigskip
\section{Proof of the Main Theorem}
\bigskip
The key point of the proof consists in showing 
that the 1-parameter family of non-linear operators, which defines
the lifts of stationary discs for the pairs $(J^t, \Gamma^t)$ of Proposition \ref{deformation}
satisfies the hypothesis of the general Implicit Function Theorem  
(see e.g. 
 \cite{KZ}, p.353)
at the disc $\hat f_{a, \lambda}$, defined in (3.8) and (3.9), and with the differential 
problem determined by  
the osculating pair $(J^0, \Gamma^0)$. After that, our main result 
will follow immediately.\par
 \medskip
So, for some fixed $0< \delta <1$, let us   consider
 the Frechet derivative  $\F'|_{\hat f_{a, \lambda}} = (\F_1|_{\hat f_{a, \lambda}} ; \F_2|_{\hat f_{a, \lambda}} )$   of the non-linear operator 
 \bigskip
 $$\F = (\F_1; \F_2): \CC^{1,\delta}(\bar \Delta, \C^{2n}) \longrightarrow \CC^\delta(\bar \Delta, \C^{2n})\times \CC^\delta(\partial \Delta, \R\times \C^{n-1} \times i\R)$$
 \par \bigskip
 \noindent
 defined by the l.h.s. of the p.d.e. system  (3.4) and of the boundary conditions
 (3.5) and (3.7). A straightforward computation shows that for any 
$\hat h = (h; k): \bar \Delta \to \C^{2n}$, one has 
$\F_1'|_{\hat f_{a, \lambda}} (\hat h) = (\H(h;k); \K(h;k))$  where
\bigskip
$$\H^\alpha(h;k)  = \frac{\partial h^\alpha}{\partial \bar \zeta}Ê\ ,\qquad 
\H^n(h;k) = \frac{\partial h^n}{\partial \bar \zeta}Ê -  \frac{i \bar a}{2} A_{\bar 1 \bar \alpha} \overline{h^\alpha}
+ \frac{i \bar a}{2} A_{\bar 1\bar \alpha} \overline{\frac{\partial h^\alpha}{\partial \zeta}} \bar \zeta\ ,$$
\bigskip
$$\K^\alpha(h;k) =  \frac{\partial k_\alpha}{\partial \bar \zeta} 
+  \frac{i a}{2} \overline{A_ {\bar \alpha \bar 1}} \zeta \overline{\frac{\partial k_n}{\partial \zeta}}
+  \frac{i \lambda}{2} \overline{A_ {\bar \alpha \bar \beta}} h^\beta\ ,\qquad 
\K^n(h;k)\=  \frac{\partial k^n}{\partial \bar \zeta} \ .$$
\par 
\bigskip
With similar computations one gets also  the components of the 
second part of the Frechet derivative  $\F_2'|_{\hat f_{a, \lambda}} (\hat h)
= (\M^0(h;k), \M^1(h;k), \dots, \M^{n-1}(h;k), \M^n(h;k)) $, i.e.  the maps 
$$\M^0(h;k): \partial \Delta 
\to \R\ , \qquad \M^\alpha(h;k): \partial \Delta 
\to \C\ ,\qquad \M^n(h;k): \partial \Delta 
\to i \R\ ,$$
with 
\bigskip
$$\M^0(h;k)(\zeta) = h^n(\zeta) + \overline{h^n(\zeta)} -  \bar a h^1(\zeta) \bar \zeta - a \overline{h^1(\zeta)} \zeta\ ,$$
\bigskip
$$\M^\alpha(h;k)(\zeta) = k_\alpha(\zeta) 
 +  \delta_{\alpha} ^1\bar a\bar \zeta k_n(\zeta) + 
 \frac{i a}{2} \overline{A_{\bar \alpha \bar 1}} \zeta \left(k_n(\zeta) - \bar k_n(\zeta)\right) +
  \phantom{aaaaaaaaaaaaaaaa}Ê$$
 $$ \phantom{aaaaaaaaaaaaaaaaaaaaaaaaa}Ê+  \lambda \left(\overline{h^\alpha(\zeta)} +
\frac{i}{2} \overline{A_{\bar \alpha\bar \beta}} h^\beta(\zeta)\right) \zeta 
-   \frac{i \lambda}{2} \overline{A_{\bar \alpha\bar \beta}} h^\beta(\zeta)\bar  \zeta \ ,$$
\bigskip
$$\M^n(h;k)(\zeta) = \bar \zeta k_n(\zeta) - \zeta \overline{k_n(\zeta)}\ .$$
\par
\bigskip
The following lemma is the key point of our proof.\par
\medskip
\begin{lem} \label{surjectivity}ÊThe linear operator $\F'|_{\hat f_{a, \lambda}}$ is surjective 
and the map  \bigskip
$$ (h,k) \in \ker \F'|_{\hat f_{a, \lambda}} \longrightarrow
\left( h^i(0)\ \text{with}\ 1 \leq i\leq n\  , 
\frac{\partial h^i}{\partial \zeta}(0)\ \text{with}\ i \geq 2, \right.\phantom{aaaaaaaaaaa}$$
$$\phantom{aaaaaaaaaaaaaaaaaaaaaaaaaaa}\left. \im \left(
\bar a \frac{\partial h^1}{\partial \zeta}(0)\right), \re\left (\frac{\partial k_n}{\partial \zeta}(0)\right)\right) \in \C^{2n-1}\times \R^2$$
\par \bigskip
is a linear isomorphism between 
$\ker \F'|_{\hat f_{a, \lambda}}$ and $\C^{2n-1}\times \R^2$.
\end{lem}
\begin{pf}  Consider the classical Cauchy-Green transform, 
i.e. the operator $T_{CG}: \CC^\delta(\Delta, \C)
\to  \CC^{1,\delta}(\Delta, \C)$, defined by
$$T_{CG}(\varphi)(\zeta) = \frac{1}{2 \pi i} \int\int_\Delta \frac{\varphi(\eta)}{\eta -\zeta} d\eta \wedge
d\bar \eta\ .$$
It is well known that it  is inverse to the $\bar \partial$-operator, i.e. satisfies $\frac{\partial T_{CG}(\varphi)}{\partial 
\zeta} = \varphi$ (see e.g. \cite{Ve}, \S I.8). With the help of $T_{CG}$, the operators $\H$ and $\K$ can be written  as 
\bigskip
$$\H^\alpha(h;k) = \frac{\partial h^\alpha}{\partial \bar \zeta}\ ,\qquad 
\H^n(h;k) = \frac{\partial}{\partial \bar \zeta}Ê\left(h^n + \frac{i \bar a}{2} A_{\bar 1\bar \alpha} \overline{h^\alpha} \bar \zeta
-  i \bar aA_{\bar 1 \bar \alpha} T_{CG}( \overline{h^\alpha})
\right)\ ,
\eqno(4.1)$$
\bigskip
$$\K^{\alpha}(h;k) = \frac{\partial}{\partial \bar \zeta}\left(
 k_\alpha +  \frac{i a}{2} \overline{A_ {\bar \alpha \bar 1}} \zeta\overline{k_n}
 +  \frac{i \lambda}{2} \overline{A_ {\bar \alpha \bar \beta}} T_{CG}( h^\beta)
 \right)  \ ,\qquad 
\K^n(h;k) = \frac{\partial k^n}{\partial \bar \zeta}\ .\eqno(4.2)$$
\par \bigskip
\noindent
From this expression and using the Cauchy-Green transform, 
one can directly shows that $\F_1'|_{\hat f_{a, \lambda}}$ 
is surjective. By linearity of the operator, this implies that  the surjectivity of $\F'|_{\hat f_{a, \lambda}}$
is proved if we can show that, for any $\Phi \in \CC^\delta(\partial \Delta, \R\times \C^{2n-1} \times i\R)$, 
there exists a solution $\hat h = (h, k)$  to the Riemann-Hilbert problem 
$$\F_1'|_{\hat f_{a, \lambda}}(\hat h) = 0\ ,\qquad \F_2'|_{\hat f_{a, \lambda}}(\hat h) = \Phi
\ .\eqno(4.3)$$ 
It follows immediately from (4.1) and (4.2) that $\hat h = (h,k)$ satisfies 
$\F_1'|_{\hat f_{a, \lambda}}(\hat h) = 0$ if and only if the components $h^\alpha$ and $k^n$ are 
holomorphic, while the components $h^n$ and $k^\alpha$ are of the form 
$$h^n =  - \frac{i \bar a}{2} A_{\bar 1\bar \alpha} \overline{h^\alpha} \bar \zeta
+  i \bar aA_{\bar 1 \bar \alpha}   \overline{\I(h^\alpha)} + \tilde h^n\ , \eqno(4.4)$$
\bigskip
$$k^\alpha =  - 
\frac{i a}{2} \overline{A_ {\bar \alpha \bar 1}} \zeta\overline{k_n}
 -  \frac{i \lambda}{2} \overline{A_ {\bar \alpha \bar \beta}} h^\beta \bar \zeta + \tilde k^\alpha
 \eqno(4.5)$$
 \par
 \bigskip
 \noindent
for some holomorphic function $\tilde h^n$ and $\tilde k^\alpha$. In (4.4) we used the symbol 
``$\I$" to denote the operator which associates to 
any holomorphic function $\phi$ on the unit disc, the unique holomorphic function $\I(\phi)$ 
such that 
$\frac{\partial \I(\phi)}{\partial \zeta} = \phi$ and $ \I(\phi)(0) = 0$. \par
From (4.4) and (4.5), it follows that (4.3) admits  a solution for an arbitrary  $\Phi$ if and only if 
there exists a holomorphic disc $( h^\alpha, \tilde h^n; \tilde k^\alpha, k^n): \bar \Delta \to \C^{2n}$ which 
satisfies the following conditions for any $\zeta \in \partial \Delta$: 
\bigskip
$$\tilde h^n(\zeta) + \overline{\tilde h^n(\zeta)}
- \frac{i \bar a}{2} A_{\bar 1 \bar \alpha} \overline{h^\alpha(\zeta)} \bar \zeta
+  \frac{i a}{2} \overline{A_{\bar 1 \bar \alpha}} h^\alpha(\zeta)  \zeta + 
i \bar a A_{\bar 1 \bar \alpha}\overline{\I(h^\alpha)(\zeta)} - 
i a \overline{A_{\bar 1 \bar \alpha}} \I(h^\alpha)(\zeta) - $$
$$ -
\bar a\ h^1(\zeta) \bar \zeta -  a\ \overline{h^1(\zeta)} \zeta = \Phi^0(\zeta)\ , \eqno(4.6)$$
\bigskip
$$\tilde k_\alpha(\zeta)  - i a \overline{A_{\bar \alpha \bar 1}} \overline{k_n(\zeta)}  \zeta
+ \delta^1_\alpha \bar a k_n(\zeta) \bar \zeta 
-  i \lambda  \overline{A_{\bar \alpha \bar \beta}} h^\beta(\zeta) \bar \zeta +
\frac{i a}{2} \overline{A_{\bar \alpha \bar 1}} \zeta k_n(\zeta) +$$
$$ + \lambda \left(\overline{h^\alpha(\zeta)} + 
\frac{ i}{2}   \overline{A_{\bar \alpha \bar \beta}} h^\beta(\zeta)\right) \zeta
= \Phi^\alpha(\zeta)\ ,
\eqno(4.7)$$
\bigskip
$$\bar \zeta k_n(\zeta) - \zeta \overline{k_n(\zeta)} =  i  \Phi^n(\zeta)\ .\eqno(4.8)$$\par
\bigskip
\noindent
We recall that, by the regularity assumed on $\Phi$,  
 each map $\Phi^i$ can be written as sum of a unique Fourier power series 
$\Phi^i(\zeta) = \sum_{m \in \Z} \Phi^i_m \zeta ^m$ (see e.g. \cite{KK}, Rmk. VII.7.5). So, if
we consider the expressions of the  holomorphic functions as 
sum of power series 
\bigskip
$$h^\alpha(\zeta) = \sum_{m \geq 0}Ê h^\alpha_m \zeta^m\ ,\quad
 \tilde h^n = \sum_{m \geq 0}Ê h^n_m \zeta^m\ ,\quad
 \tilde k_\alpha= \sum_{m \geq 0}Ê k_{\alpha,m} \zeta^m\ ,\quad k^n(\zeta) 
 =  \sum_{m \geq 0}Ê k_{n,m} \zeta^m\ ,$$
 \par
 \bigskip
 \noindent
we obtain the following conditions on the coefficients $h^i_m$ and $k_{i,m}$, with  $m \geq 0$
(for reader's convenience, we point out that the equations (4.9) are obtained from (4.8), equations (4.10) - (4.13) come from (4.7) and 
(4.14) and (4.15) are consequences of  (4.6)):
\bigskip
$$ k_{n,m} = i \Phi^n_{m-1} \ \ \text{for all} \  m \geq 3\ ,\ \  k_{n,2} - \overline{k_{n,0}} = i\Phi^n_{1}\ ,\ \ 
k_{n,1} - \overline{k_{n,1}} =   i \Phi^n_{0}\ ,\eqno(4.9)$$
\bigskip
$$\overline{h^\alpha_{m}} = \frac{1}{\lambda} \Phi^\alpha_{1-m} 
+\frac{ i a}{\lambda} \overline{A_{\bar \alpha \bar 1}} \overline{k_{n,m}}
\ \text{for all} \  m \geq 3 \ ,\ \ 
\overline{h^\alpha_2} - i \overline{A_{\bar \alpha \bar \beta}} h^\beta_0 = 
\frac{i a}{\lambda} \overline{A_{\bar \alpha, \bar 1}} \overline{k_{n,2}} 
- \frac{\delta^1_\alpha}{\lambda} \bar a k_{n,0}  + \frac{1}{\lambda} \Phi^\alpha_{-1}
\eqno(4.10) $$
\bigskip
$$k_{\alpha, 0} =   i a \overline{A_{\bar \alpha \bar 1}} \overline{k_{n,1}}
- \delta^1_\alpha \bar a k_{n,1} +
i \lambda  \overline{A_{\bar \alpha \bar \beta}} h^\beta_1 
- \lambda \overline{h^\alpha_1} + \Phi^\alpha_0\ ,\eqno(4.11)$$
\bigskip
$$k_{\alpha, 1} =   i a \overline{A_{\bar \alpha \bar 1}} \overline{k_{n,0}}
- \delta^1_\alpha \bar a k_{n,2}  -  \frac{i a}{2} \overline{A_{\bar \alpha \bar 1}} k_{n,0}   
 + i \lambdaÊ\overline{A_{\bar \alpha \bar \beta}}
\left( h^\beta_2 - \frac{h^\beta_0}{2} \right) - \lambda \overline{h^\alpha_0} 
+ \Phi^\alpha_1\ ,\eqno(4.12)$$
\bigskip
$$k_{\alpha, m} =   
- \delta^1_\alpha \bar a k_{n,m+1} -  \frac{i a}{2} \overline{A_{\bar \alpha \bar 1}} k_{n,m-1}   +
i \lambda  \overline{A_{\bar \alpha \bar \beta}}\left( h^\beta_{m+1} -
\frac{h^\beta_{m-1}}{2}\right)
+ \Phi^\alpha_m\ \text{for all} \  m \geq 2\ ,\eqno(4.13)$$ 
\bigskip
$$ h^n_m =  i a \overline{A_{\bar 1 \bar \alpha}} h^\alpha_{m-1}\left(
 \frac{1}{m} - \frac{1}{2}\right) + \bar a h^1_{m+1} + \Phi^0_m\ \ \ \ \text{for all} \  m \geq 2\ ,\eqno(4.14)$$
 \bigskip
$$h^n_1 = \frac{i a}{2} \overline{A_{\bar 1 \bar \alpha}} h^\alpha_{0}
 + \bar a h^1_2 + a \overline{h^1_0} + \Phi^0_1\ ,\ \ 
h^n_0 + \overline{h^n_0} = \bar a h^1_1 + a \overline{h^1_1} + \Phi^0_0\ .\eqno(4.15)$$
\par \bigskip
\noindent
This system is immediately seen to be solvable by  substitutions 
for any choice of the Fourier coefficients 
$\Phi^i_m$ and this concludes the proof of the surjectivity. \par
\medskip
Let us now consider $\ker \F'|_{\hat f_{a, \lambda}}$. It is clear that it is linearly  isomorphic with the space of 
holomorphic discs $(h^\alpha, \tilde h^n, \tilde k^\alpha,  k^n): \bar \Delta \to \C^{2n}$, 
whose power series coefficients satisfy the above equations with all terms $\Phi^i_m$'s
set equal to $0$. In particular, it can be checked 
 that any solution is uniquely determined by the 
values of  $h^i_0$,  for $1\leq i \leq n$, by the real numbers 
$\re(k_{n,1})$ and $\im(\bar a h^1_1)$ and  by the complex 
numbers $h^i_1$ for $2 \leq i \leq n$ (to check this claim, recall 
that  $A_{\bar 1\bar 1} = 0$ and  then solve all the
equations by substitutions in the following order: $(4.15)$, $(4.10)_2$ for $\alpha = 1$, 
$(4.9)$, $(4.10)$ and then all the others). From this, 
the second claim follows immediately.
\end{pf}
\bigskip
From Lemma \ref{surjectivity}, we see that the continuous family of operators
$$(\F^{(t)}, \V^{(t)}): \CC^{1,\delta}(\bar \Delta) \longrightarrow
\CC^\delta(\bar \Delta, \C^{2n})\times \CC^\delta(\partial \Delta,  \C^{n-1} \times \R) \times 
\C^{n}\times \R^{2n-1}\times \R$$
given by the differential operators $\F^{(t)}$,  defining the lifts of stationary 
discs for the pairs $(J^t, \Gamma^t)$, and by the evaluation maps
$$ \hat f = (f,g) \overset{\V^{(t)}}\longrightarrow \left( f^i(0),\ \ \ 
\frac{a^{-1}\cdot \frac{\partial f^i}{\partial \zeta}(0)}{ 
\re  \left(a^{-1} \cdot \frac{\partial f^1}{\partial \zeta}(0)\right)}\ \ \ , \re\left(
\frac{\partial g_n}{\partial \zeta}(0)\right)\right)$$
has invertible Frechet derivative at any disc $\hat f_{a, \lambda}$. \par
\medskip
By the general Implicit
Function Theorem, there exists a stationary disc passing through any point $z_o$ on the inward  real normal to $\Gamma^0$ 
in a sufficiently small neighborhood of the origin and tangent to any vector 
of  an open neighborhood of 
 $\span_\C\left\{ \left. \frac{\partial}{\partial z^\alpha}\right|_z\ \right\}_{1 \leq \alpha \leq n-1\ }  $
for any pair $(J^t, \Gamma^t)$, $0 \leq t < \epsilon$, for $\epsilon$ sufficiently small (moreover, this obviously remains true for any point 
$z_o$ in an open  convex cone with vertex at the origin containing the real inward normal of $\Gamma$ if this cone contains no real lines). 
By Proposition \ref{deformation}, our main theorem follows.\par

\bigskip
\bigskip
\font\smallsmc = cmcsc8
\font\smalltt = cmtt8
\font\smallit = cmti8
\hbox{\parindent=0pt\parskip=0pt
\vbox{\baselineskip 9.5 pt \hsize=3.1truein
\obeylines
{\smallsmc
Andrea Spiro
Dip. Matematica e Informatica
Universit\`a di Camerino
Via Madonna delle Carceri
I-62032 Camerino (Macerata)
ITALY
}\medskip
{\smallit E-mail}\/: {\smalltt andrea.spiro@unicam.it}
}
\hskip 0.0truecm
\vbox{\baselineskip 9.5 pt  \hsize=3.7truein
\obeylines

{\smallsmc
Alexandre Sukhov
U.S.T.L.
Cit\'e Scientifique
59655 Villeneuve d'Asq Cedex
FRANCE
}\medskip
{\smallit E-mail}\/: {\smalltt  sukhov@agat.univ-lille1.fr}
}
}

\end{document}